\documentclass{amsart}
\usepackage{xypic, amssymb, amsfonts, amsbsy, amsthm, amsmath, amscd, latexsym, stmaryrd, epic, eepic,eucal}

\DeclareMathAlphabet{\mathpzc}{OT1}{pzc}{m}{it}

\newenvironment{dem}{\begin{proof}[\bf Proof]}{\end{proof}}
\newtheorem{Thm}{Theorem}[section]
\newtheorem{theorem}{\bf Theorem}[section]
\newtheorem{lemma}[theorem]{\bf Lemma}
\newtheorem{propos}[theorem]{\bf Proposition}

\newtheorem{claim}[theorem]{\bf Claim}

\theoremstyle{definition}
\newtheorem{defi}[theorem]{\bf Definition}
\newtheorem{rmk}[theorem]{\bf Remark}

\newtheorem{exm}[theorem]{\bf Example}

\newcommand{\A}{\mathbb A}
\newcommand{\B}{\text{B}}

\newcommand{\M}{\mathcal M}

\newcommand{\Mcal}{\mathcal M}
\newcommand{\N}{\mathbb N}
\newcommand{\Ncal}{\mathcal N}
\newcommand{\Of}{\mathcal O}
\newcommand{\Pro}{\mathbb P}

\newcommand{\Z}{\mathbb Z}
\newcommand{\aut}{\text{Aut}}

\newcommand{\cart}{\ar @{} [dr] |{\Box}}

\newcommand{\Hom}{\text{Hom}}

\newcommand{\id}{\text{id}}

\newcommand{\Obj}{\text{Obj}}

\newcommand{\spe}{\text{Spec}}

\newcommand{\ov}{\overline}
\newcommand{\ord}{\text{ord}}

\newcommand{\becomes} {\xymatrix @C=12pt @R=4pt  {
& &  &\\
& \ar[r] & &\\
& & &
}}

\newcommand{\opentwo} {\xymatrix @C=20pt @R=20pt  {
{\bullet} \text{\tiny 2}
}}

\newcommand{\dtwoone} {\xymatrix @C=20pt @R=20pt  {
*=0{\; \bullet \text{\tiny 1}} \ar@(ul,dl)@{-}
}}

\newcommand{\dtwotwo} {\xymatrix @C=20pt @R=20pt  {
*=0{ \text{\tiny 1} \bullet} \ar@{-}[r] & *=0{\bullet \text{\tiny 1}}
}}

\newcommand{\ctwoone} {\xymatrix @C=20pt @R=20pt  {
*=0{\bullet} \ar@(ul,dl)@{-} \ar@(ur,dr)@{-}
}}

\newcommand{\ctwotwo} {\xymatrix @C=20pt @R=20pt  {
*=0{\bullet} \ar@{-}[r] \ar@(ul,dl)@{-} & *=0{\bullet \text{\tiny 1}}
}}

\newcommand{\ptwoone} {\xymatrix @C=20pt @R=20pt  {
*=0{\bullet} \ar@{-}[r] \ar@(ul,dl)@{-} & *=0{\bullet}  \ar@(ur,dr)@{-}
}}

\newcommand{\ptwotwo} {\xymatrix @C=20pt @R=20pt  {
*=0{\bullet} \ar@{-}[r] \ar@{-} @/_/[r]  \ar@{-} @/^/[r]  & *=0{\bullet} 
}}

\newcommand{\pthreeone} {\xymatrix @C=20pt @R=20pt  {
*=0{\bullet} \ar@{-}[r] \ar@{-}[d] \ar@{-}[dr] & *=0{\bullet} \ar@{-}[dl] \ar@{-}[d] \\
*=0{\bullet} \ar@{-}[r] & *=0{\bullet}
}}

\newcommand{\pthreetwo} {\xymatrix @C=20pt @R=20pt  {
*=0{\bullet} \ar@{-} @/^/[d] \ar@{-} @/_/[d] \ar@{-}[r] & *=0{\bullet} \ar@{-} @/^/[d] \ar@{-} @/_/[d] \\
*=0{\bullet} \ar@{-}[r] & *=0{\bullet}
}}

\newcommand{\pthreethree} {\xymatrix @C=12pt @R=4pt  {
*=0{\bullet} \ar@{-} @/^/[dd] \ar@{-} @/_/[dd] \ar@{-}[dr] & & \\
& *=0{\bullet} \ar@{-}[r] & *=0{\bullet} \ar@(ur,dr)@{-}\\
*=0{\bullet} \ar@{-}[ur] & &
}}

\newcommand{\pthreefour} {\xymatrix @C=20pt @R=4pt  { *=0{\bullet} \ar@(ul,dl)@{-} \ar@{-}[r]& *=0{\bullet} \ar@{-} @/^/[r] \ar@{-} @/_/[r]& *=0{\bullet} \ar@{-}[r] & *=0{\bullet} \ar@(ur,dr)@{-}}}

\newcommand{\pthreefive} {\xymatrix @C=9pt @R=3pt { & *=0{\bullet} \ar@(ur,ul)@{-} \ar@{-}[dd] & \\ & &  \\ & *=0{\bullet} \ar@{-}[dl] \ar@{-}[dr] & \\ *=0{\bullet} \ar@(l,d)@{-} & & *=0{\bullet} \ar@(r,d)@{-} }}

\newcommand{\cthreeone} {\xymatrix @C=4pt @R=8pt  {
 & *=0{\bullet} \ar@{-}[dl] \ar@{-} @/_/[dl] \ar@{-} @/^/[dr] \ar@{-} [dr]& \\
 *=0{\bullet} \ar@{-}[rr] & & *=0{\bullet}
}}

\newcommand{\cthreetwo} {\xymatrix @C=12pt @R=4pt  {
*=0{\bullet} \ar@{-} @/^/[dd] \ar@{-} @/_/[dd] \ar@{-}[dr] &  \\
& *=0{\bullet} \ar@(ur,dr)@{-}\\
*=0{\bullet} \ar@{-}[ur] &
}}

\newcommand{\cthreethree} {\xymatrix @C=20pt @R=4pt  { *=0{\bullet} \ar@(ul,dl)@{-} \ar@{-} @/^/[r] \ar@{-} @/_/[r]& *=0{\bullet} \ar@{-}[r] & *=0{\bullet} \ar@(ur,dr)@{-}}}

\newcommand{\cthreefour} {\xymatrix @C=9pt @R=3pt { & *=0{\bullet} \ar@(ur,ul)@{-} \ar@{-}[dl] \ar@{-}[dr] & \\ *=0{\bullet} \ar@(l,d)@{-} & & *=0{\bullet} \ar@(r,d)@{-} }}

\newcommand{\cthreefive} {\xymatrix @C=20pt @R=4pt  { & & \\
*=0{\bullet} \ar@{-}[r] \ar@{-} @/^/[r] \ar@{-} @/_/[r]& *=0{\bullet} \ar@{-}[r] & *=0{\bullet} \ar@(ur,dr)@{-}\\
& &
}}

\newcommand{\cthreesix} {\xymatrix @C=12pt @R=4pt  {
*=0{\bullet} \ar@{-} @/^/[dd] \ar@{-} @/_/[dd] \ar@{-}[dr] & & \\
& *=0{\bullet} \ar@{-}[r] & *=0{\bullet \text{\tiny 1}} \\
*=0{\bullet} \ar@{-}[ur] & &
}}

\newcommand{\cthreeseven} {\xymatrix @C=20pt @R=4pt  { *=0{\bullet} \ar@(ul,dl)@{-} \ar@{-}[r]& *=0{\bullet} \ar@{-} @/^/[r] \ar@{-} @/_/[r]& *=0{\bullet} \ar@{-}[r] & *=0{\bullet \text{\tiny 1}}}}

\newcommand{\cthreeeight} {\xymatrix @C=9pt @R=3pt { & *=0{\stackrel{1}{\bullet}} \ar@{-}[dd] & \\ & &  \\ & *=0{\bullet} \ar@{-}[dl] \ar@{-}[dr] & \\ *=0{\bullet} \ar@(l,d)@{-} & & *=0{\bullet} \ar@(r,d)@{-} }}

\newcommand{\cfivespecial}{\xymatrix @C=20pt @R=20pt {
& *=0{\bullet} \ar@{-}[dl] \ar@{-}[d] \ar@{-}[dr] & \\
*=0{\bullet} \ar@{-}[r] \ar@{-}[d] & *=0{\bullet} \ar@{-}[r] \ar@{-}[d] & *=0{\bullet} \ar@{-}[d] \\
*=0{\bullet} \ar@{-} @/^/[r] \ar@{-} @/_/[r] & *=0{\bullet} & *=0{\bullet} \ar@(dr,dl)@{-}
}}

\newcommand{\pfivespecialone}{\xymatrix @C=16pt @R=5pt {
*=0{\bullet} \ar@{-} @/^/[dd] \ar@{-} @/_/[dd] \ar@{-}[r]  & *=0{\bullet} \ar@{-}[r] \ar@{-}[dd] & *=0{\bullet} \ar@{-}[dr] \ar@{-}[dd] & &\\
& & & *=0{\bullet} \ar@{-}[r] & *=0{\bullet} \ar@(ur,dr)@{-}\\
*=0{\bullet} \ar@{-}[r] & *=0{\bullet} \ar@{-}[r] & *=0{\bullet} \ar@{-}[ur]& &
}}

\newcommand{\pfivespecialtwo}{\xymatrix @C=22pt @R=22pt {
*=0{\bullet} \ar@{-}[r] \ar@{-}[d] \ar@{-}[dr]& *=0{\bullet} \ar@{-}[dl] \ar@{-}[dr] & \\
*=0{\bullet} \ar@{-}[d] & *=0{\bullet} \ar@{-}[r] \ar@{-}[d] & *=0{\bullet} \ar@{-}[d] \\
*=0{\bullet} \ar@{-} @/^/[r] \ar@{-} @/_/[r] & *=0{\bullet} & *=0{\bullet} \ar@(dr,dl)@{-}
}}

\newcommand{\pfivespecialthree}{\xymatrix @C=20pt @R=20pt {
& *=0{\bullet} \ar@{-}[dl] \ar@{-}[r] \ar@{-}[dr] & *={\bullet} \ar@{-}[ddl]\\
*=0{\bullet} \ar@{-}[r] \ar@{-}[d] & *=0{\bullet} \ar@{-}[r] \ar@{-}[ur] & *=0{\bullet} \ar@{-}[d] \\
*=0{\bullet} \ar@{-} @/^/[r] \ar@{-} @/_/[r] & *=0{\bullet} & *=0{\bullet} \ar@(dr,dl)@{-}
}}

\newcommand{\lszeroone}{\xymatrix @C=8pt @R=8pt  {  *=0{\bullet} \ar@{-}[dr] & &  *=0{\bullet} \ar@{-}[dl]\\ &  *=0{\bullet} \ar@{-}[d] & \\  *=0{\bullet} \ar@{-}[r] &  *=0{\bullet} \ar@{-}[r] &  *=0{\bullet}
}}

\newcommand{\lszerotwo}{\xymatrix @C=8pt @R=8pt  {  *=0{\bullet} \ar@{-}[dr] & &  *=0{\bullet} \ar@{-}[dl]\\ *=0{\bullet} \ar@{-}[r] &  *=0{\bullet} \ar@{-}[r] &  *=0{\bullet}
}}

\newcommand{\lszerothree}{\xymatrix @C=8pt @R=8pt  { & *=0{\bullet} \ar@{-}[d] &  *=0{\bullet} \ar@{-}[d]& \\ *=0{\bullet} \ar@{-}[r] &  *=0{\bullet} \ar@{-}[r] & *=0{\bullet} \ar@{-}[r] &  *=0{\bullet}
}}

\newcommand{\lsone}{\xymatrix @C=20pt @R=4pt  { *=0{\bullet} \ar@{-} @/^/[r]^{1} \ar@{-} @/_/[r]_{2}& *=0{\bullet} \ar@{-} @/^/[r]^3 \ar@{-} @/_/[r]_4 & *=0{\bullet}
}}

\newcommand{\lsonepone}{\xymatrix @C=20pt @R=4pt  { *=0{\bullet} \ar@{-} @/^/[r]^{1} \ar@{-} @/_/[r]_{2} & *=0{\bullet} \ar@{-}[r]^p & *=0{\bullet} \ar@{-} @/^/[r]^3 \ar@{-} @/_/[r]_4 & *=0{\bullet}
}}

\newcommand{\lsoneptwo}{\xymatrix @C=20pt @R=4pt {& *=0{\bullet} \ar@{-}[dr]^3 \ar@{-}[dd]^p & \\ *=0{\bullet} \ar@{-}[ur]^{1} \ar@{-}[dr]_2 & & *=0{\bullet}\\ & *=0{\bullet} \ar@{-}[ur]_4&
}}

\newcommand{\lstwo} {\xymatrix @C=12pt @R=4pt  {
*=0{\bullet}  \ar@{-}[dr]^{1} &  & \\
& *=0{\bullet} \ar@{-} @/^/[r]^3 \ar@{-} @/_/[r]_4 & *=0{\bullet} \\
*=0{\bullet} \ar@{-}[ur]_{2} & &
}}

\newcommand{\lsthree} {\xymatrix @C=30pt @R=16pt  { *=0{\bullet} \ar@{-}[r]^{2} \ar@{-} @/^/[r]^{1} \ar@{-} @/_/[r]_{3}& *=0{\bullet} \ar@{-}[r]^{4} & *=0{\bullet} }}

\newcommand{\lsthreep} {\xymatrix @C=12pt @R=4pt  {
*=0{\bullet} \ar@{-} @/^/[dd] \ar@{-} @/_/[dd] \ar@{-}[dr]^{p} & & \\
& *=0{\bullet} \ar@{-}[r] & *=0{\bullet} \\
*=0{\bullet} \ar@{-}[ur] & &
}}

\newcommand{\lsfour}{\xymatrix @C=16pt @R=16pt  {
& *=0{\bullet}  \ar@{-}[d]^{1} & \\
*=0{\bullet}  \ar@{-}[r]^{2} & *=0{\bullet}  \ar@{-}[r]_{4} &*=0{\bullet} \\
& *=0{\bullet}  \ar@{-}[u]^{3} &
}}

\newcommand{\graphzfour} {\xymatrix @C=15pt @R=15pt  {
& & & & *=0{\bullet} \ar@{-}[d] \ar@(ur,ul)@{-} & & \\ 
&*=0{\bullet} \ar@{-}[r] \ar@{-}[d] \ar@{-}[ddrr] & *=0{\bullet} \ar@{-} @/^/[r] \ar@{-} @/_/[r]& *=0{\bullet} \ar@{-}[r] & *=0{\bullet} \ar@{-}[r] & *=0{\bullet} \ar@{-}[ddll] \ar@{-}[d] & \\
*=0{\bullet} \ar@{-}[r] \ar@(ul,dl)@{-} &*=0{\bullet} \ar@{-}[d]  & & & &*=0{\bullet}  \ar@{-} @/^/[d] \ar@{-} @/_/[d] & \\
&*=0{\bullet}  \ar@{-} @/^/[d] \ar@{-} @/_/[d]  & & *=0{\bullet} & & *=0{\bullet} \ar@{-}[d] & \\
&*=0{\bullet}  \ar@{-}[d]  & & & & *=0{\bullet} \ar@{-}[d] & *=0{\bullet} \ar@{-}[l] \ar@(ur,dr)@{-} \\
&*=0{\bullet} \ar@{-}[r]  \ar@{-}[uurr] & *=0{\bullet} \ar@{-}[r]  & *=0{\bullet} \ar@{-} @/^/[r] \ar@{-} @/_/[r] & *=0{\bullet} \ar@{-}[r] & *=0{\bullet}  \ar@{-}[uull] & \\
& & *=0{\bullet} \ar@{-}[u]  \ar@(dr,dl)@{-} & & & & \\ 
}}

\newcommand{\graphonethreethree} {\xymatrix @C=15pt @R=15pt  {
&*=0{\bullet} \ar@{-}[r] \ar@{-}[d] \ar@{-}[dddr] & *=0{\bullet} \ar@{-} @/^/[r] \ar@{-} @/_/[r]& *=0{\bullet} \ar@{-}[dddrr]&  & & \\
*=0{\bullet} \ar@{-}[r] \ar@(ul,dl)@{-} &*=0{\bullet} \ar@{-}[d]  & & & & & \\
&*=0{\bullet}  \ar@{-} @/^/[d] \ar@{-} @/_/[d]  & & & & & \\
&*=0{\bullet}  \ar@{-}[d]  & *=0{\bullet} & *=0{\bullet}  \ar@{-}[l] \ar@(ur,dr)@{-} & & *=0{\bullet} \ar@{-}[d] & *=0{\bullet} \ar@{-}[l] \ar@(ur,dr)@{-} \\
&*=0{\bullet} \ar@{-}[r]  \ar@{-}[ur] & *=0{\bullet} \ar@{-}[r]  & *=0{\bullet} \ar@{-} @/^/[r] \ar@{-} @/_/[r] & *=0{\bullet} \ar@{-}[r] & *=0{\bullet}  \ar@{-}[ulll] & \\
& & *=0{\bullet} \ar@{-}[u]  \ar@(dr,dl)@{-} & & & & \\ 
}}

\newcommand{\maxloop} {\xymatrix @C=25pt @R=15pt  {
					& *=0{\bullet} \ar@{-}[d] \ar@(ur,ul)@{-}	& *=0{\bullet} \ar@{-}[d] \ar@(ur,ul)@{-}	& *=0{\bullet} \ar@{-}[d] \ar@(ur,ul)@{-}	& & *=0{\bullet} \ar@{-}[d] \ar@(ur,ul)@{-} & *=0{\bullet} \ar@{-}[d]  \ar@(ur,ul)@{-} & \\
\ar@(ul,dl)@{-} *=0{\bullet} \ar@{-}[r]_1 	& *=0{\bullet} \ar@{-}[r]_2 	& *=0{\bullet} \ar@{-}[r]_3 	& *=0{\bullet} \ar@{.}[rr] 	& & *=0{\bullet} \ar@{-}[r]_{g-2}  & *=0{\bullet} \ar@{-}[r]_{g-1} & *=0{\bullet} \ar@(ur,dr)@{-}
}}

\input xy
\xyoption{all}

\begin{document}

\title{Moduli Stacks of Curves with a Fixed Dual Graph}

\author{Dan Edidin}
\author{Damiano Fulghesu}
\address{Department of Mathematics, University of Missouri, Columbia, MO 65211}
\address{IRMA -- Universit\'e de Strasbourg et CNRS 7, rue R. Descartes
67084 Strasbourg, Cedex; France}
\email{edidind@missouri.edu, fulghesu@math.u-strasbg.fr}
\thanks{The first author was partially supported by NSA grant H98230-08-1-0059.}
\begin{abstract}
  Let $\M_g(\Gamma)$ be the stack of stable curves of genus $g$ with a
  given dual graph $\Gamma$ and let $\ov{\M}_g(\Gamma)$ be its closure
  in $\ov{\M}_g$. We give an explicit desingularization of
  $\ov{\M}_g(\Gamma)$ and we study the one-dimensional substack of
  $\ov{\M}_g$ of curves with at least $3g - 4$ nodes.
\end{abstract}

\maketitle

\section{Introduction}

In his famous paper \cite{mumenum}, Mumford introduced a
stratification of the stack $\ov{\M}_g$
given by the number of nodes. The stratum, $\mathcal{M}^n_g$,
corresponding to curves with $n$ nodes has pure dimension $3g-3
-n$,  but
is not irreducible. The irreducible components of 
$\M_g^n$ are indexed by weighted stable graphs (see Definition
\ref{stable}). 
Precisely, if $\Gamma$ is weighted stable graph with $n$ edges and
(weighted) genus $g$ then substack $\M_g(\Gamma)$ parametrizing curves
with dual graph $\Gamma$ is an irreducible component of
$\M_g^n$.
The substacks $\M_g(\Gamma)$ and
their closures in $\ov{\M}_g$ give a combinatorial decomposition of
the stack $\ov{\M}_g$, and a natural question is to describe the
irreducible components of Mumford's stratification.
Unfortunately, the combinatorics is rather complicated due to the rapid
growth of the number of possible graphs $\Gamma$ for increasing genus
and number of nodes.

In this paper we give an explicit desingularization of
$\ov{\M}_g(\Gamma)$ in the category of Deligne-Mumford stacks (Theorem \ref{normalization}).
This allows to verify whether $\ov{\M}_g(\Gamma)$ is singular, knowing the
automorphism group of $\Gamma$ (see Remark \ref{des}).

We then focus on the 1-stratum of $\ov{\M}_g$, that is the stratum
corresponding to curves with $3g-4$ nodes. We first prove in Section
\ref{toprmk} that the 1-stratum is connected in $\ov{\M}_g$.
(Theorem \ref{connectedness}). As a corollary we conclude that the
strata $\ov{\M}_g^n$ are connected for all $n$. 

In Section \ref{1strat} we classify the irreducible components
$\ov{\M}_g(\Gamma)$ of the 1-stratum of $\ov{\mathcal M}_g$ (for a
different approach to the same problem see \cite{Zintl}). These
components are all Deligne-Mumford stacks whose coarse moduli space is
$\mathbb P^1$ (see Proposition \ref{coarse} and Remark
\ref{caseh1}). Moreover, we describe their normalizations as gerbes
over an orbifold (Theorems \ref{structurenormh0} and
\ref{structurenormh1}). Although we do not enumerate the number of
irreducible components of $\ov{{\mathcal M}}_g^{(3g-4)}$, we show that there
are only 5 possible residual orbifolds (Remark \ref{residualorbifolds}
and Theorem \ref{structurenormh1}). We also give a local presentation
for the singular points of $\ov{\M}_g(\Gamma)$ (Proof of Proposition
\ref{coarse}), and give examples to show that all possible residual
orbifolds are realized.

\section{Preliminaries}

\subsection{Deligne-Mumford stacks} We will work with  stacks over
a noetherian base scheme $S$. This means in particular that a stack
$\mathcal X$ will be considered equipped with a morphism
$$ \psi: \mathcal X \to S.
$$
Stacks are defined as categories fibered in groupoids over a site
(with some extra conditions). The base scheme $S$ represents category of schemes of finite presentation over $S$ equipped with the \'etale topology.
For basic definitions of stacks we refer to \cite{DM}, \cite{Art}, \cite{LMB} and Appendix in \cite{Vis}. Here we gather together a few basic facts.

\begin{defi}
A morphism $F: \mathcal X \to \mathcal Y$ between two stacks over $S$
is {\em representable} if for every scheme $T$ and for every morphism
$T \to \mathcal Y$, the fibered product $\mathcal X \times_{\mathcal
  Y} T  \to T$ is a scheme.
\end{defi}

Many concepts about morphisms of schemes may be applied to representable morphisms of stacks.

\begin{defi}\label{morphisms}
  Let $\mathbf{P}$ be a property of morphisms of schemes that is
  stable under base change and of local nature on the target
  (e.g. flat, smooth, \'etale, surjective, unramified, normal, locally
  of finite type, locally of finite presentation). Then we say that a
  representable morphism of stacks $\mathcal X \to \mathcal Y$ has
  property $\mathbf{P}$ if for every morphism $T \to \mathcal Y$, the
  morphism of schemes deduced by base change $g: T \times_{\mathcal Y}
  \mathcal X$ has that property.
\end{defi}

\begin{defi} \label{def.DMstack}
A stack ${\mathcal X}$ is a {\em Deligne-Mumford (DM)  stack} if the
following conditions are satisfied.

1) The diagonal $\Delta \colon {\mathcal X} \to {\mathcal X} \times_S
{\mathcal X}$ is representable, quasi-compact and separated.

2) There exists an \'etale surjective morphism $U \to {\mathcal X}$
where $U$ is a scheme.

The scheme $U$ is called and {\em atlas} for ${\mathcal X}$.

An {\em algebraic space} is a DM stack which is equivalent to a sheaf.
\end{defi}
\begin{rmk}
The representability of the diagonal implies that any morphism $U \to
{\mathcal X}$ with $U$ a scheme is representable and the morphism $U
\to {\mathcal X}$ of Definition \ref{def.DMstack} is \'etale and
surjective in the sense of Definition \ref{morphisms}.
\end{rmk}

\begin{rmk}
If $f \colon {\mathcal X} \to {\mathcal Y}$ is a morphism of DM stacks
then
$f$ has property $\mathbf{P}$ if for some (and hence every) \'etale
atlas
$U \to {\mathcal X}$ the morphism of schemes $U \times_{\mathcal Y}
{\mathcal X} \to U$ has property ${\mathbf P}$.
\end{rmk}

Let $\mathbf{P}$ be a property of schemes, local in the \'etale topology (e.g. regular, smooth, normal, reduced), then we say that a Deligne-Mumford stack has property $\mathbf{P}$ if and only if the atlas $U$ satisfies $\mathbf{P}$.

\begin{rmk}\label{allmorphisms}
The structure morphism $\psi: \mathcal X \to S$ is not representable unless $\mathcal X$ is a scheme. So, according to the given definition, we cannot say that $\mathcal X$ satisfies $\mathbf{P}$ if and only if $\psi$ satisfies $\mathbf{P}$. However if $\mathbf{P}$ is a property of local nature, at source and target, for the \'etale topology (e.g. flat, smooth, \'etale, unramified, normal, locally of finite type, locally of finite presentation), then we can extend the definitions for morphisms of DM stacks  which are not necessarily representable (see p.100 \cite{DM}).
\end{rmk}

\begin{defi}
A stack ${\mathcal X}$ is separated over $S$ if the diagonal $\Delta
\colon {\mathcal X} \to {\mathcal X} \times_S {\mathcal X}$ is a
finite representable morphism.
\end{defi}

By \cite[Theorem 2.7]{EHKV} every Deligne-Mumford stack ${\mathcal X}$ admits a finite
surjective morphism $Z \to {\mathcal X}$ with $Z$ a scheme. Using this
fact we can define the notions of proper and finite morphisms of DM
stacks. Our definition of proper morphism is equivalent to the one
given in \cite{DM}.
\begin{defi}
A morphism of DM stacks ${\mathcal X} \to {\mathcal Y}$ is proper (resp. finite) if for some (and
hence all) finite surjective morphism $Z \to {\mathcal X}$ with $Z$ a
scheme, the composite morphism $Z \to {\mathcal X} \to {\mathcal Y}$
is a representable proper (resp. finite) morphism.
\end{defi}

\begin{propos} \label{prop.finiteisrep}
Let $F \colon {\mathcal X} \to {\mathcal Y}$ be finite surjective morphism of DM
stacks which is faithful (i.e. the functor $F$ is a faithful functor)
then $F$ is representable.
\end{propos}

To prove Proposition \ref{prop.finiteisrep} we begin with a Lemma.

\begin{lemma}\label{rep}
Let $F: \mathcal X \to \mathcal Y$ be a morphism between two DM
stacks. Then $F$ is weakly representable\footnote{A morphism is {\em
    weakly representable} if for every morphism $T \to {\mathcal Y}$
  with $T$ an algebraic space, the fiber product is represented by an
  algebraic space.} if and only if $F$ is faithful.
\end{lemma}

\begin{dem}
From Corollary (8.1.2) \cite{LMB} we have that $F$ is weakly representable if and only if the diagonal morphism
\begin{eqnarray*}
\Delta_F : \mathcal X &\to& \mathcal X \times_{\mathcal Y} \mathcal X \\
x \in \Obj(\mathcal X) &\mapsto&  (x,x, \id_{F(x)})\\
(f: x \to x') &\mapsto& (f,f,F(f))
\end{eqnarray*}
is fully faithful. This condition means exactly that two morphisms
$$
f,g \in \Hom_{\mathcal X}(x,x')
$$
are equal if and only if $F(f)=F(g)$.
\end{dem}

\begin{dem}{Proof of Proposition \ref{prop.finiteisrep}.}
By Lemma \ref{rep} the map faithful functor $F \colon {\mathcal X} \to
{\mathcal Y}$ is weakly representable. Let $T \to {\mathcal X}$ be a
morphism from a scheme and let ${\mathcal X}_T$ denote the fiber
product $T \times_{\mathcal Y} {\mathcal X}$. Since $F$ is weakly
representable we know that ${\mathcal X}_T$ is an algebraic space. To
prove that $F$ is representable we need to show that ${\mathcal X}_T$
is actually a scheme. Working locally on $T$ we may assume that $T$ is
affine. Let $Z \to {\mathcal X}$ be a finite surjective morphism from
a scheme. Since $Z \to {\mathcal Y}$ is finite, surjective and
representable, the fiber product $Z_T = T \times_{\mathcal Y} Z$ is
represented by a scheme and the map $Z_T \to T$ is finite and
surjective. Since $T$ is assumed affine, the scheme $Z_T$ is also
affine. Now the morphism $Z_T \to {\mathcal X}_T$ is, by base change, a
finite surjective morphism of algebraic spaces. Chevalley's theorem
for algebraic spaces \cite[Chapter III, Theorem 4.1]{Knu} implies that ${\mathcal X}_T$ is an
affine scheme as well.
\end{dem}

\begin{defi} Let $\mathcal X$ be an DM stack. A geometric point of $\mathcal X$ is a morphism
$$
\spe (K) \xrightarrow{x} \mathcal X
$$
where $K$ is an algebraically closed field. From any such map we can deduce an object $\xi$ in $\mathcal X (\spe (K))$. Let $G_x$ be the automorphism group of $\xi$, we have a monomorphism
$$
\B G_x \xrightarrow{\text{rg}_x} \mathcal X
$$
We call $G_x$ the {\it stabilizer} of $x$ and $\text{rg}_x$ the {\it residual gerbe} of $x$. If $G_x$ is not trivial, we say that $x$ is a {\it stacky point}\footnote{In literature it is also called {\it twisted point}.} of $\mathcal X$.
\end{defi}

Let $x: \spe ( K ) \to \mathcal X$ be a geometric point of a Deligne-Mumford stack $\mathcal X$, and let $U_{x}$ be an \'etale scheme-theoretic neighborhood of $x$. We have a lifting of $x$ to $U_{x}$, so we define
$$
\widehat{\mathcal O}_{x, \mathcal X} := \widehat{\mathcal O}_{x, U_{x}}
$$

We recall the following fundamental Lemma (see \cite{AV00}, 2.2.2 and 2.2.3). It states that, locally in the \'etale topology, every separated Deligne-Mumford stack is a quotient stack by a finite group.

\begin{lemma}\label{localquotient}
Let $\mathcal X$ be a separated Deligne-Mumford stack, and $X$ its coarse moduli space. There is an \'etale covering $\{ X_{\alpha} \to X\}_{\alpha \in I}$, such that for each $\alpha \in I$ there is a a scheme $U_{\alpha}$ and a finite group $G_{\alpha}$, acting on $U_{\alpha}$, with the property that the pullback $\mathcal X \times_X X_{\alpha}$ is isomorphic to the quotient stack $[U_{\alpha}/G_{\alpha}]$. Moreover, if $X$ is noetherian, $X_{\alpha}$ is the coarse moduli space of $[U_{\alpha}/G_{\alpha}]$ (therefore it is isomorphic to the geometric quotient $U_{\alpha}/G_{\alpha}$).
\end{lemma}

\begin{rmk}
All the stacks we will consider are proper, therefore separated.
\end{rmk}

\begin{propos}\label{localinvariant}
Let $\mathcal X$ be a separated Deligne-Mumford stack with a noetherian coarse moduli space $X$. Let $x: \spe (K) \to \mathcal X$ be a geometric point and $G_x$ its stabilizer. Then

(a) there is a natural action of $G_x$ on $\widehat{\mathcal O}_{x, \mathcal X}$;

(b) $\widehat{\mathcal O}_{x, X} \cong \left( \widehat{\mathcal O}_{x, \mathcal X} \right)^{G_x}$.
\end{propos}

\begin{dem}

(a) Following the notation of Lemma \ref{localquotient}, we consider an \'etale neighborhood $[U_{\alpha}/G_{\alpha}]$ of $x$ in $\mathcal X$. $G_x$ is a subgroup of $G_{\alpha}$ that stabilizes a lifting of $x$ to $U_{\alpha}$, therefore we have an induced action of $G_{x}$ on $\widehat{\mathcal O}_{x, \mathcal X} = \widehat{\mathcal O}_{x, U_{\alpha}}$.

(b) We have that an \'etale neighborhood of $x$ in $X$ is the geometric quotient $U_{\alpha}/G_{\alpha}$. We can also choose $U_{\alpha}=\spe (R)$ for some ring $R$. A groupoid which represents $\left[ \spe (R) / G \right]$ is
\begin{equation*}
W= \xymatrix{
\left( \spe (R) \times _{k} G \right) \ar@<1ex>[rr]^{\pi_1} \ar@<-1ex>[rr]_{\gamma}& & \spe (R)
}
\end{equation*}
where $\pi_1$ is the first projection and $\gamma$ is the action. Now, from \cite{KM} Proposition 5.1, we have that the geometric quotient is $\spe \left( R^W \right)$. In our case $R^W=R^G$, so the coarse moduli space of $\left[ \spe (R) / G \right]$ is $\spe \left( R^G \right)$. In order to conclude we must show that the natural morphism
$$
\spe \left( R^{G_x} \right) \to \spe \left(R^G \right)
$$
is \'etale on $x$, and this is given by \cite{SGA} Expos\'e V, Proposition 2.2.
\end{dem}

\subsection{Quotients of DM  stacks}

In the following we will consider quotients of DM stacks by an action of a finite group acting on it.
More precisely by {\it group} we mean a sheaf in groups over the base category. Usually the base is the category of schemes over a base scheme $S$ that we simply call $S$ (in this case the base category is the stack structure of $S$). The main reference is \cite{Rom}. Here we recall some basic definitions (see (loc. cit.) Definitions 2.1 and 2.3).

\begin{defi}
Let $\M$ be a stack over a base scheme $S$, and let $G$ be a sheaf in
groups over $S$. Let $m$ be the multiplication of $G$, and $e$ its
unit section. An {\it action} of $G$ on $\M$ is a morphism of stacks
$\mu: G \times \M \to \M$ with strictly commutative diagrams
$$
\xymatrix{
G \times G \times \M \ar[r]^{m \times \text{id}} \ar[d]_{id \times \mu} & G \times \M \ar[d]^{\mu} \\
G \times \M \ar[r]^{\mu} & \M
}
\qquad
\xymatrix{
G \times \M \ar[r]^{\mu} & \M \\
\M \ar[u]^{e \times \text{id}} \ar[ur]_{\text{id}}&
}
$$
We say that $\M$ is a $G$-stack. 
\end{defi}

Any stack over $S$ can be seen as a $G$-stack over $S$ through the trivial action.

\begin{defi}
Let $\M$ and $\mathcal{N}$ be two $G$-stacks and $\psi: \mathcal N \to \M$ a morphism. We say that $\psi$ is a $G$-morphism ($\psi \in \hom_{G\text{-stacks}}(\mathcal N, \M)$), if the diagram
 $$
\xymatrix{
G \times  \mathcal N \ar[r]^{\mu_{\mathcal{N}} } \ar[d]_{id \times \psi} &  \mathcal N \ar[d]^{\psi} \\
G \times \M \ar[r]^{\mu_{\M}} & \M
}
$$
is strictly commutative.
\end{defi}

\begin{rmk}
If we consider groupoids over $S$ instead of stacks the above diagrams are 2-diagrams satisfying some ``higher associativity" condition (see \cite{Rom} Definition 1.3). In the case of stacks we require that the action is strict. 
\end{rmk}

\begin{defi} Let $G$ be a sheaf in groups over $S$ and let $\M$ be a $G$-stack over $S$. A {\it quotient stack} $\M / G$ is a stack that 2-represents the 2-functor
$$
F(\mathcal N) = \hom_{G\text{-stacks}}(\mathcal N, \M).
$$
\end{defi}

\begin{propos} (\cite{Rom} Theorem 3.3)
Let $G$ be a sheaf in groups over $S$. and $\M$ a $G$-stack over $S$. Then there exists a quotient stack $\M / G$ and its formation commutes  with base change on $S$.
\end{propos}

Under some hypothesis on $G$ we can extend the above proposition to DM stacks (see Theorem 4.1 \cite{Rom}).

\begin{theorem}\label{quotgroup}
  Let $G$ be an \'etale group scheme over $S$. Let $\M$ be a DM $G$-stack over $S$. Then the quotient stack $\M / G$ is a DM stack.
\end{theorem}

\section{Topological classes of stable curves}

\begin{defi}\label{stable}
We call {\it weighted graph}, a graph with a natural number on each vertex. Given a weighted graph, we call $V(\Gamma)$ the set of vertices, $E(\Gamma)$ the set of edges and $w: V(\Gamma) \to \N$ the assignment of weights. Let $c=|V(\Gamma)|$, $n=|E(\Gamma)|$ and $h= \sum_{v \in V(\Gamma)} w(v)$. For each vertex $v$ in $|V(\Gamma)|$ we call {\it degree} of $v$ ($\text{deg}(v)$) the number of edges starting from $v$ (loops counting twice). We also define {\it multiplicity} of $v$ to be the number $\text{mult}(v)=3w(v)+ \text{deg}(v)$. We call {\it weighted genus} of a weighted graph the number \footnote{This is not the standard genus of a graph that we can find in literature, more precisely we are adding the total weight $h$.}
\begin{equation}\label{graphrel}
g= h + n - c + 1.
\end{equation}
We say that a weighted graph $\Gamma$ is {\it stable} if $\Gamma$ is connected, its genus is $\geq 2$ and for every $v \in V(\Gamma)$ we have $\text{mult}(v) \geq 3$.

We call loop graph a graph whose cycles are only loops.
\end{defi}

\begin{rmk}
There is a subtlety that should be considered at this point. We are
think of graphs as unlabeled because 
we want to consider a unique topological class of a stable curve that does not depend on labeling. However in the definition of graph and in the following definition of automorphism, we need to label vertices and edges by fixing the sets $V(\Gamma)$ and $E(\Gamma)$.
\end{rmk}

\begin{defi}
An automorphism of a weighted graph $\Gamma$ is the following set of data:
\begin{enumerate}
\item a one-to-one correspondence $f: V(\Gamma) \to V(\Gamma)$ such that $w(f(v))=w(v)$ for all vertices of $\Gamma$;

\item a one-to-one correspondence $g: E(\Gamma) \to E(\Gamma)$ such that for all $e \in E(\Gamma)$ the two vertices connected by $g(e)$ are images of the vertices connected by $e$;

\item an element of $\Z^{l(\Gamma)}_2$ where $l(\Gamma)$ is the number of the loops of $\Gamma$ \footnote{We are thinking that we can flip each loop and in general that there are two ways to send a loop into another one}.
 \end{enumerate}
We call $\aut(\Gamma)$ the group of automorphisms of $\Gamma$.
\end{defi}

Given a stable curve $C$ of weighted genus $g$ we can associate a weighted stable graph $\Gamma$ of genus $g$ by setting \footnote{Notice that a node belongs to at most two components and we set a loop when a component has self intersection.}
\begin{eqnarray*}
V(\Gamma)&=&\{ \text{ irreducible components of $C$ } \}, \\
E(\Gamma)&=&\{ \text{ nodes of $C$ } \} \\
& \text{and}&\\
w: V(\Gamma) &\to& \N \\
w(v) &:=& \text{ genus of the normalization of $v$}
\end{eqnarray*}

We call $\Gamma$ the {\it topological class} of $C$.

As an example we give now the table that we get  for genus $g=2$. Each graph corresponds to a topological class of a curve of genus 2. We put in columns graphs with a given total weight $h$ and in lines we fix the number of components.

\begin{center}

\begin{tabular}{|l|c|c|c|}

\hline
$g=2$ & $h=2$ & $h=1$ & $h=0$\\

\hline

 $c=1$ & $\opentwo$ & $\qquad$ $\dtwoone$ & \begin{tabular}{c}
                                                                $\;$\\
                                                                $\ctwoone$\\
                                                                $\;$
                                                        \end{tabular}\\
\hline

$c=2$ & $\dtwotwo$ & $\qquad$ $\ctwotwo$ & \begin{tabular}{c}
                                                        $\;$\\
                                                        $\qquad$ $\ptwoone$ $\qquad$\\
                                                        $\;$\\
                                                        $\ptwotwo$\\
                                                        $\;$
                                                        \end{tabular}\\

\hline

\end{tabular}

\end{center}

In general the total weight runs from 0 to $g$ while the number of components is at most $2g-2$. We add here also the case $g=3$ but we just write the number of stable weighted graphs in each square.

\begin{center}
\begin{tabular}{|l|c|c|c|c|}
\hline
$g=2$ & $h=3$ & $h=2$ & $h=1$ & $h=0$\\
\hline
$c=1$ & 1 & 1 & 1 & 1 \\
\hline
$c=2$ & 1 & 3 & 4 & 4 \\
\hline
$c=3$ & 1 & 3 & 6 & 5 \\
\hline
$c=4$ & 1 & 2 & 3 & 5 \\
\hline
\end{tabular}
\end{center}

\begin{rmk}
Instead of the total weight we can also use for labeling the columns the standard genus $g_s$ of the graph, since $h+g_s=g$.
\end{rmk}

\begin{defi}
For any integer $g \geq 2$ and any weighted stable graph $\Gamma$ of genus $g$ we call $\M_g(\Gamma)$ the full (locally closed) substack of $\ov{\M}_g$ consisting of stable curves with topological class $\Gamma$.

We also call $i$-stratum of $\ov{\M}_g$ the union of the closures in $\ov{\M}_g$ of $\M_g(\Gamma)$ such that $\dim \left( \M_g(\Gamma) \right)=i$. This is the stratum of $\ov{\M}_g$ consisting of curves with $3g-3 - i$ nodes
\end{defi}

\begin{rmk} Given a stable weighted graph $\Gamma$ we have
$$
\dim ( \M_g(\Gamma) ) = \sum_{v \in V(\Gamma)} (\text{mult}(v) - 3).
$$
Notice also that in the above tables the anti-diagonals preserve the number of edges (i.e. of nodes), therefore graphs in the same anti-diagonal correspond to components of the same dimension.
\end{rmk}

\begin{propos}\label{smooth}
For every stable graph $\Gamma$ of genus $g$ the stack $\M_g(\Gamma)$
is smooth over $\spe (\Z)$.
\end{propos}

\begin{dem}
We know by \cite{DM} that $\ov{\M}_g$ is a smooth stack over
$\spe (\Z)$. To prove that
the substack $\M_g(\Gamma)$ is smooth it suffices to show that the substack
$\M_g(\Gamma)$ is formally smooth.
Let $x \colon \spe (k) \to \M_g(\Gamma)$ be a geometric point corresponding
to a stable curve $C_x$ of topological class $\Gamma$. 

The complete local ring to $\ov{\M}_g$ at $x$ is the complete local ring of
the universal deformation space ${\mathcal M}$ of the curve $C_x$.
By \cite{DM} this ring is ${\mathfrak o}_k[[t_1,\ldots , t_{3g-3}]]$. Here
${\mathfrak o}_k = k$ if the characteristic is 0
and ${\mathfrak o}_k$ is the unique complete
regular local ring with residue field $k$ and maximal ideal $p{\mathfrak o}_k$
if the characteristic is $p$.

Moreover, we may choose the $t_i$'s such that if ${\mathcal C} \to
{\mathcal M}$ is the universal curve then complete local ring of
${\mathcal C}$ at the nodes of $C_x$ is isomorphic to $${\mathfrak
o}_k[[u_i,v_i,t_1, \ldots , t_{3g-3}]]/(u_iv_i -t_i).$$

The complete local ring of $\M_g(\Gamma)$ at $x$ is the quotient of ${\mathcal
O}_{x,\ov{\M}_g}$ by the ideal corresponding to deformations that
preserve the nodes. From the description of the complete local rings
to ${\mathcal C}$ at the nodes of $C$ we see that this ideal is $(t_1,
\ldots ,t_r)$ where $r$ is the number of edges of $\Gamma$. Hence
$\M_g(\Gamma)$ is smooth.
\end{dem}

\section{Normalization of the substacks $\ov{\M}_g(\Gamma)$}

Given a vertex $v$ in a graph $\Gamma$, we call $\widehat{E}(v)$ the (ordered) set of edges meeting $v$ and considering loops twice.

\begin{defi} Given a stable weighted graph $\Gamma$ of weighted genus $g$ we define the stack
$$
\Ncal_g(\Gamma) := \left( \prod_{v \in V(\Gamma)} \Mcal_{w(v), \widehat{E}(v)} \right)
$$
and the natural 1-morphism
$$
\Ncal_g(\Gamma) \xrightarrow{\pi_{\Gamma}} \Mcal_g(\Gamma)
$$
induced by gluing sections corresponding to the same edge.

Moreover we define $\ov{\Ncal}_g(\Gamma)$ as the stack
$$
\ov{\Ncal}_g(\Gamma):= \left( \prod_{v \in V(\Gamma)} \ov{\Mcal}_{w(v), \widehat{E}(v)} \right)
$$
and extend $\pi_{\Gamma}$:
$$
\ov{\Ncal}_g(\Gamma) \xrightarrow{\ov{\pi}_{\Gamma}} \ov{\Mcal}_g(\Gamma).
$$
\end{defi}

\begin{propos}\label{properties}
The 1-morphisms $\pi_{\Gamma}$ and $\ov{\pi}_{\Gamma}$ are (a)
representable; (b) finite; (c) unramified; (d) surjective. Moreover (e) the 1-morphism $\pi_{\Gamma}$ is \`etale.
\end{propos}

\begin{dem} 
  The morphism $\ov{\pi}_\Gamma$ is a composition of clutching
  morphisms in the sense of Knudsen \cite{Knud2}. Knudsen proved that
  clutching morphisms are representable, finite and unramified
  \cite[Corollary 3.9]{Knud2}. To prove (d) we have to check that for
  all curves $C \to T$ in $\M_g(\Gamma)$ there exists an \'etale
  covering $T' \to T$ such that $C_{T'} \to T'$ (obtained by base
  change) is isomorphic to the image of some object in $\mathcal N
  _g(\Gamma)$.

  First of all fix a geometric point $\spe (k) \to T$. We can choose a
  smooth point on each irreducible component of the fiber $C_k$
  defining sections $\{ s_{v} \}_{v \in V(\Gamma)}$. We know that
  there exists an \'etale covering $T_1 \to T$ where the $s_{v}$
  extend. So we get a curve $C_1 \to T_1$ of topological class
  $\Gamma$ where we have labeled irreducible components with
  $V(\Gamma)$.

  Let us now consider the normalization $\widehat{C_1} \to C_1 \to
  T_1$. The pre-image of the relative singular locus of $C_1$ defines
  a divisor $\widehat{D_1}\subset \widehat{C_1}$ and an \'etale
  covering $\widehat{D_1} \to T_1$ of degree $\sum_{v V(\in \Gamma)}
  |\widehat{E}(v)|$. If this covering is trivial then we can choose
  $T':=T_1$ and $C_{T'}=C_1$ \footnote{The sections $\{ s_v\}_{v \in
      V(\Gamma)}$ rigidify the components, but not necessarily the
    nodes.}.

  Otherwise let $H_1, H_2, \dots, H_r$ be all the irreducible
  components of $\widehat{D_1}$ such that each morphism $H_i \to T_1$
  is not trivial. For every $1 \leq i \leq r$, let $q_1$ be the degree
  of $H_i \to T_1$. We call $q=\sum_{i=1}^r q_i$ the excess covering
  number of $\pi_1 \to T_1$. Let
$$
T_2= \coprod_{i=1}^r H_i
$$
and consider the cartesian diagram
$$
\xymatrix{
\widehat{C_2} \ar[r] \ar[d] \cart &  \widehat{C_1} \ar[d]\\
C_2 \ar[r] \ar[d]_{\pi_2} \cart &  C_1 \ar[d]^{\pi_1}\\
T_2 \ar[r]^{\sigma_1} & T_1
}
$$

The morphism $\widehat{C_2} \to T_2$ admits at least the identity section, therefore $\pi_2: C_2 \to T_2$ has a strictly smaller excess covering number.

After a finite number of steps we get the required $T'$ and $C_{T'}$.

(e) We have that $\mathcal N_g(\Gamma)$ is smooth because it is the
product of smooth stacks, moreover we have proved in proposition
\ref{smooth} that $\M_g(\Gamma)$ is smooth, therefore $\pi_{\Gamma}$
is a finite representable surjective morphism between a smooth
(hence Cohen-Macaulay)
stack and a smooth (hence regular) stack of the same dimension.
Hence it is flat
\cite[cf. Remark 3.11]{liu}.
From Proposition 17.6.1 \cite{egaq} we have that a flat
and unramified morphism of relative dimension 0 is necessarily
\'etale.
\end{dem}

There is a natural action of $\aut(\Gamma)$ on $\Ncal_g(\Gamma)$ (that extends to  $\ov{\Ncal}_g(\Gamma)$) which is consistent with the action on $\Gamma$ so the quotient stacks
$$
\Ncal_g(\Gamma)/\aut(\Gamma) \qquad \ov{\Ncal}_g(\Gamma)/\aut(\Gamma)
$$
are well defined from Theorem \ref{quotgroup}.

Let us consider a geometric point of $\M_g(\Gamma)$, that is to say a
stable curve $C$ over an algebraically closed field $K$ of topological class $\Gamma$. We have the
following diagram
\begin{equation*}
\xymatrix{
\spe (K) \times_{\Mcal_g(\Gamma)} \Ncal_g(\Gamma) \cart \ar[r] \ar[d] & \Ncal_g(\Gamma) \ar[d]^{\pi_{\Gamma}} \\
\spe (K) \ar[r] & \M_g(\Gamma)
}
\end{equation*}

\begin{propos}\label{exun}
There is a natural isomorphism
$$
\spe (K) \times_{\Mcal_g(\Gamma)} \Ncal_g(\Gamma) \cong \aut(\Gamma)
$$
\end{propos}

\begin{dem}
 It is enough to check the isomorphism over $\spe (K)$. Let $C$ be the
 curve over $K$ with topological class $\Gamma$ defined by $\spe (K) \to
 \Mcal_g(\Gamma)$. The objects of the groupoid
$$ \Ncal_g (\Gamma)_K (\spe (K) ) = \left( \spe (K)
\times_{\Mcal_g(\Gamma)} \Ncal_g(\Gamma) \right)(\spe (K))
$$ are pairs $(\widehat C , \alpha)$ where $\widehat C$ is an object
in $\Ncal_g (\Gamma)(\spe (K))$ and $\alpha$ is an isomorphism between
$C$ and $\pi_{\Gamma}(\widehat C)$.

The isomorphisms between two objects $(\widehat C , \alpha)$ and
$(\widehat C' , \alpha')$, are isomorphisms $g: \widehat C \to
\widehat C'$ such that the diagram
\begin{eqnarray*}
\xymatrix{
\pi_{\Gamma}(\widehat C) \ar[rr]^{\pi_{\Gamma}(g)}& & \pi_{\Gamma}(\widehat C')\\
& C \ar[ul]^{\alpha} \ar[ur]_{\alpha'}&
}
\end{eqnarray*}
commutes. That is to say $\pi_{\Gamma}(g) = \alpha' \alpha^{-1}$. By
representability of $\pi_{\Gamma}$ we have at most one isomorphism $g$
having this property. In particular this means, as we already know from Proposition \ref{properties} (c), that
$\mathcal{N}_g(\Gamma)_K$ is a set of points.

Let us now fix an object $(\widehat C , \alpha)$ and take $\gamma \in \aut(\Gamma)$ different from the identity. Let us write
$$
\widehat C = \coprod_{v \in V(\Gamma)} C_v
$$
where $C_v$ is a curve in $\Mcal_{w(v), \widehat{E}(v))}$. Let us define
$$
\gamma(\widehat C) := \coprod_{v \in V(\Gamma)} C_v
$$
where now $C_v$ is a curve in $\Mcal_{w(\gamma(v)), \gamma(\widehat{E}(v))}$.

There is an isomorphism $\beta$  between $\pi_{\Gamma}(\widehat{C})$ and $\pi_{\Gamma}(\gamma(\widehat{C}))$ sending each component $v$ to $\gamma(v)$ and each node $e$ to $\gamma(e)$. Define $\gamma(\alpha)= \beta \alpha: C \to \pi_{\Gamma}(\gamma(\widehat{C}))$. Since there are no isomorphisms in $\mathcal N _g (\Gamma)$ between $\widehat{C}$ and $\gamma(\widehat{C})$ whose image through $\pi_{\Gamma}$ is $\beta$, we have that $(\gamma(\widehat{C}), \gamma(\alpha))$ defines a different point in $\mathcal N_g(\Gamma)_K$.

On the other hand let $(\widehat{C}_2, \alpha_2)$ be another object in $\mathcal N_g(\Gamma)_K$ such that $\alpha_2 \alpha^{-1}$ cannot be lifted to an isomorphism $\widehat{\beta}: \widehat{C}_2 \to \widehat{C}$.

Let us now consider the following diagram
\begin{equation*}
\xymatrix{
\widehat{C} \ar[rr]^{\delta} \ar[d] & & \widehat{C}_2 \ar[d]\\
\pi_{\Gamma}(\widehat{C}) \ar[rr]^{\alpha_2 \alpha^{-1}} & & \pi_{\Gamma}(\widehat{C}_2)\\
& C \ar[ur]_{\alpha_2} \ar[ul] ^{\alpha} &
}
\end{equation*}
where vertical maps are normalization morphisms which are consistent with the labeling of $\Gamma$.  There exists a unique isomorphism $\delta$ making the diagram commute. Moreover, since $\delta$ cannot be an isomorphism in $\Ncal _g(\Gamma)_K$, it cannot preserve all nodes and components. Therefore $\delta$ defines an element $\gamma$ in $\aut(\Gamma)$ different from identity.

So we can conclude that $(C_2, \alpha_2)$ is isomorphic to $(\gamma(\widehat{C}), \gamma(\alpha))$.
\end{dem}

\begin{propos}
There is an isomorphism
$$
\left[ \Ncal_g(\Gamma)/ \aut(\Gamma) \right] \cong \Mcal_g(\Gamma)
$$
\end{propos}

\begin{dem}
From the proof of Proposition \ref{exun} we get that the induced morphism
$$
\left[ \Ncal_g(\Gamma)/ \aut(\Gamma) \right] \to \Mcal_g(\Gamma)
$$
is faithful, that is to say representable (from Lemma \ref{rep}). From Proposition \ref{properties} we have that it is an \'etale covering and using Proposition \ref{exun} we get it is of degree one. Hence it is an isomorphism.
\end{dem}

We still have a morphism
$$
\varphi_{\Gamma}: [\ov{\Ncal}_g(\Gamma)/ \aut(\Gamma)] \to \ov{\M}_g(\Gamma)
$$
but it is in general far from being an isomorphism. The main reason is that we could get extra automorphisms at the points in the closure.

However, $\ov{\Ncal}_g(\Gamma)$ is normal, since it is the product of normal stacks, so also the quotient stack is normal. This means that $\varphi_{\Gamma}$ factorizes through the normalization of $\ov{\M}_g(\Gamma)$.

\begin{theorem}\label{normalization}
The morphism
$$
\varphi_{\Gamma}: [\ov{\Ncal}_g(\Gamma)/ \aut(\Gamma)] \to \ov{\M}_g(\Gamma)
$$
is a normalization for $\ov{\M}_g(\Gamma)$.
\end{theorem}

\begin{dem}
The map $\varphi_\Gamma$ is faithful on automorphism groups by
Proposition \ref{exun}. It is also finite and surjective so by
by Proposition \ref{prop.finiteisrep} it is representable.  Since
normalization commutes with \'etale base change \cite[Proposition
18.12.15]{egaq} we may assume that $\varphi_\Gamma$ is a map of
schemes. Since the source is normal and the map is finite and generically an isomorphism, $\varphi_\Gamma$ must be the normalization by the
universal property of normalizations.
\end{dem}

\begin{rmk} \label{des}
Note also that since $[\ov{\Ncal}_g/(\Gamma)/\aut(\Gamma)]$ is non-singular, the map
$\varphi_{\Gamma}$ is a desingularization of $\ov{\M}_g(\Gamma)$. Moreover $\ov{\M}_g(\Gamma)$ is smooth if and only if $\varphi_{\Gamma}$ is an isomorphism.
\end{rmk}

\section{Topological considerations on the 1-stratum}\label{toprmk}

Throughout this section all graphs will be weighted and stable.

We want to prove that the 1-stratum of genus $g \geq 2$ stable curves is connected.

First of all notice that each of the irreducible components of the 1- stratum comes from a graph $\Gamma$ having either of the following properties:
\begin{enumerate}
\item all vertices of $\Gamma$ have weight 0, exactly one of them has degree 4 and all the others have degree 3;

\item exactly one among the vertices of $\Gamma$ has weight 1 (and degree 1) and all the others have degree 3 and weight 0.
\end{enumerate}

Fix a graph $\Gamma$ coming from the 1-stratum, then the geometric points we have to add to $\ov{\Mcal}_g(\Gamma)$ correspond to the graphs obtained from $\Gamma$ in the following ways:
\begin{enumerate}
\item if we are in the first case we split the vertex of degree 4 into two vertices and add an edge between them as in the following example:

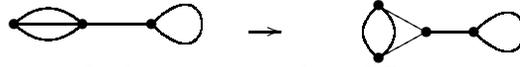
\begin{tabular}{ccc}
$\cthreefive$ & $\becomes$ & $\pthreethree$
\end{tabular}

we can do this in at most three different ways depending on the symmetries of the graph;

\item if we are in the second case we change the weight 1 into 0 and add a loop. For example:
\end{enumerate}
\bigskip

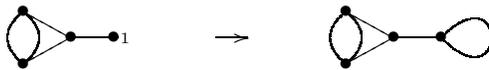
\begin{tabular}{ccc}
$\cthreesix$ & $\becomes$ & $\pthreethree$
\end{tabular}

\bigskip

\begin{defi}
We call a transformation of graphs like the two above {\it pop}. On the other hand we call {\it shrink} the inverse of pop, that is to say shrinking any edge of graph from the 0-stratum and adding a weight 1 if we shrink a loop.
\end{defi}
\begin{rmk}
It is a straightforward computation to check that pop and shrink preserve stability and genus.
\end{rmk}

\begin{Thm}\label{connectedness}
For any $g \geq 2$, the 1-stratum of $\ov{\M}_g$ is connected.
\end{Thm}

\begin{dem}

By applying a pop we see that any curve in $\ov{\Mcal}_g (\Gamma)$ contains a geometric point of the 0-stratum. In order to prove connectedness  of the 1-stratum it is enough to show that the 0-stratum is connected through components of the 1-stratum. If a graph $\Gamma$ of the 0-stratum can be turned into another graph $\Gamma'$ of the 0-stratum through a finite sequence of shrinks and pops, then $\M_{g}(\Gamma)$ is connected with $\M_g(\Gamma')$ through the 1-stratum. We do this by proving two claims.

\end{dem}

\begin{claim}
Any graph of the 0-stratum can be turned into a loop graph of the 0-stratum by applying a finite sequence of shrink-pop transformations.
\end{claim}

\begin{dem}
Take a graph $\Gamma$ of the 0-stratum and a cycle of $\Gamma$ which is not a loop. Shrink one of the edges of the cycle getting a smaller cycle. Then apply the only pop that preserves the shorter cycle and continue until we get a loop.
\end{dem}

\begin{rmk} We cannot simply take the largest cycle because in the process it is not guaranteed that in the process other cycles are not enlarged. Actually in the process we could also obtain isomorphic graphs at different steps: the relevant  part is keeping track of the cycle we are reducing in order to take out a loop.
\end{rmk}

\begin{claim}
Any loop graph of (weighted) genus $g \geq 2$ can be turned into the following  loop graph
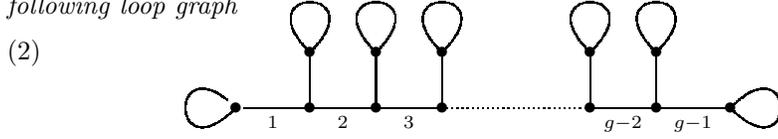
\begin{equation}\label{maxloop}
\maxloop
\end{equation}
by applying a finite sequence of shrink-pop transformations.
\end{claim}

\begin{dem}
Choose a path of maximal length of the graph. If the graph is not like (\ref{maxloop}), then we have somewhere the following subgraph (here the horizontal edges belong to the longest path)
$$
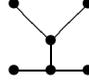
\lszeroone
$$

\bigskip

that we can transform as follows
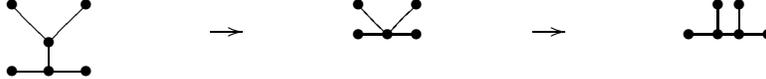
\begin{equation*}
\lszeroone \qquad \becomes \qquad \lszerotwo \qquad \becomes \qquad \lszerothree
\end{equation*}

obtaining a graph whose horizontal edges belong to the (now unique) maximal path. It is clear that after a finite numbers of shrink-pop transformations we get the graph (\ref{maxloop}).
\end{dem}

{\bf Remark} As a trivial consequence, we have that, for any $1 \leq i \leq 3g-4$, the $i$-stratum of  $\ov{\M}_g$ is connected.

\section{A description of the irreducible components of the 1-stratum}\label{1strat}

In this section we work over a base field $S= \spe (k)$. We assume that $\text{char}(k) \neq 2$ and that the field $k$ contains the cubic roots of $-1$. These are made in order to consider standard results of the action of $S_4$ on $\ov{\M}_{0,4}$.

Let $\Gamma$ be a stable weighted graph of weighted genus $g$ and $3g-3$ edges. As explained in Section \ref{toprmk}, we have two possibilities:
\begin{enumerate}
\item all vertices of $\Gamma$ have weight 0, exactly one of them has degree 4 and all the others have degree 3 (total weight $h=0$);

\item exactly one among the vertices of $\Gamma$ has weight 1 (and degree 1) and all the others have degree 3 and weight 0 (total weight $h=1$).
\end{enumerate}

\subsection{Cross ratio} Before starting with the case $h=0$ we need to point out some considerations about the cross ratio. More precisely we define a morphism

$$
\Pro^1 \backslash \{ 0, 1, \infty \} \cong  \Mcal_{0,4} \xrightarrow{cr} \Pro^1 \backslash  \{ 0, 1, \infty \}
$$

that sends $\lambda \in \Pro^1 \backslash \{ 0, 1, \infty \}$ to the cross ratio of $\{ 0, 1, \infty, \lambda \}$

This map is an isomorphism and can be extended to an isomorphism to
$$
\Pro^1 \cong  \ov{\Mcal}_{0,4} \xrightarrow{cr} \Pro^1
$$

So once we have chosen the order of the marked points of a rational curve in $\ov{\Mcal}_{0,4}$ there is a natural way to associate the cross ratio. Even if we should consider the isomorphism $cr$ we call $\lambda \in \ov{\Mcal}_{0,4}$ the cross ratio of the associated marked curve.

Let us now consider the action of $S_4$ on $\ov{\Mcal}_{0,4}$ that permutes the marked points. It is known that the orbit of $\lambda$ for the action of $S_4$ is generically of order 6:
$$
\lambda, \quad \frac{1}{\lambda}, \quad \frac{\lambda-1}{\lambda}, \quad \frac{\lambda}{\lambda-1}, \quad 1- \lambda, \quad \frac{1}{1-\lambda}
$$

and the generic stabilizer is the normal subgroup $V_4 \cong \Z/2\Z \times \Z/2\Z$ generated by $(12)(34)$ and $(13)(24)$. Therefore, we have an induced action of $S_3 \cong S_4/V_4$, that, up to isomorphisms, corresponds to the permutation of $\{ 0,1, \infty \}$ (that is to say the first three points).

On $\Pro^1$ we have the exceptional orbits:
\begin{eqnarray*}
&& o_1:= \left \{ 0,1, \infty \right \}; \\
&& o_2:= \left \{ \frac{1}{2}, 2, -1 \right \};\\
&& o_3:= \left\{ \xi,  \xi^{-1} \right \}.
\end{eqnarray*}
where $\xi \neq -1$ and $\xi^3 = -1$.
The generic stabilizer for the action of $S_3$ is the identity. The stabilizers of the two exceptional orbits $o_1$ and $o_2$ are the three subgroups generated respectively by $(1,2), (1,3), (2,3)$. Finally the stabilizer of the orbit $o_3$ is the subgroup generated by $(123)$.

\begin{rmk} Two curves with topological class $\Gamma$ are isomorphic if only if the four nodes on the same component have the same cross ratio for some order consistent with $\Gamma$.
\end{rmk}

\subsection{The case $h(\Gamma)=0$.}

Let us call $v_0$ the point with degree 4 and let $e_1, e_2, e_3, e_4$ the four edges in $\widehat{E}(\Gamma)$ (where loops are counted twice) ending in it.

Let $\sigma: \aut(\Gamma) \to S_4$ be the map that sends every automorphism $\alpha \in \aut(\Gamma)$ to the corresponding permutation $\sigma_{\alpha} \in S_4$ of the four edges. Clearly we have
\begin{eqnarray*}
&& \sigma(\alpha^{-1})=\sigma(\alpha)^{-1};\\
&& \sigma(\alpha \alpha') = \sigma(\alpha) \sigma(\alpha').
\end{eqnarray*}
Therefore, $\sigma$ is a group homomorphism.

\begin{defi}
Let us define
$$
R(\Gamma):= \aut(\Gamma)/ \sigma^{-1}(V_4) \subseteq S_3
$$
\end{defi}

\begin{theorem}\label{structurenormh0}
The normalization of $\ov{\M}_g(\Gamma)$ is a $\sigma^{-1}(V_4)$-gerbe (in the sense of \cite{BN} Proposition 4.6) over the orbifold
$$
\left [  \ov{\M}_{0,4} / R(\Gamma) \right ]
$$
\end{theorem}

\begin{dem}
We have
$$
\ov{\Ncal}_g(\Gamma)  \cong \ov{\M}_{0,4}
$$
since all the other factors are $\M_{0,3}=\spe (k)$.

From Theorem \ref{normalization} we get that the normalization of $\ov{\M}_g(\Gamma)$ is $\left [ \ov{\M}_{0,4}  / \aut(\Gamma) \right]$ where the generic stabilizer is the normal subgroup $\sigma^{-1}(V_4)$. We conclude by using the same argument of \cite{BN} Proposition 4.6.
\end{dem}

\begin{rmk}\label{residualorbifolds} The smooth Deligne-Mumford stack $\left [  \ov{\M}_{0,4} / R(\Gamma) \right ]$ is the residual orbifold of the normalization of $\ov{\M}_g(\Gamma)$. In order to define uniquely such orbifolds, it is enough to give the order of the orbifold points, since the coarse moduli space is $\Pro^1$ for all of them (see \cite{Vis} Proposition 2.11).

Now, $R(\Gamma)$ is a subgroup of  $S_3$, so we have only four possibilities: $\{ \id \}, \Z/2\Z, \Z/3\Z, S_3$. For each possibility, we write $\left [  \ov{\M}_{0,4} / R(\Gamma) \right ] = [\underline{a}|\underline{b}]$ ,where $\underline{a}$ is the set of orders of the points in the closure and $\underline{b}$ is the set of the orders of the remaining orbifold points. We summarize everything in the following chart:

\begin{tabular}{|c|c|}
\hline
$R(\Gamma)$ & $\left [  \ov{\M}_{0,4} / R(\Gamma) \right ]$\\
\hline
$\{ \id \}$ & $[1,1,1| \emptyset ]$\\
\hline
$\Z/2\Z$ & $[1,2|2]$\\
\hline
$\Z/3\Z$ & $[1|3,3]$\\
\hline
$S_3$ & $[2|2,3]$\\
\hline
\end{tabular}

\end{rmk}

\begin{propos} \label{coarse}
The coarse moduli space of $\ov{\M}_g(\Gamma)$ is $\mathbb P^1$.
\end{propos}

It may seem counterintuitive that the non-normal (and hence singular)
stacks
$\ov{\M}_g(\Gamma)$ can have smooth coarse moduli spaces. However,
as the following the example shows, it is  possible for the geometric
quotient of a non-normal scheme by a finite group to be normal.
\begin{exm}
There is an obvious $\Z_2$ action on
$\A^2$
which exchanges the coordinates. If $X = \spe k[x,y]/(xy)$ then the
action of $\Z_2$ on the invariant subscheme $X$ exchanges the two
irreducible components. The $\Z_2$-invariant subring of $k[x,y]/(xy)$
is isomorphic to the ring $k[t]$ where $t = x+ y$. Thus 
$\A^1$ is the coarse moduli space of the non-normal stack $[X/\Z_2]$.
\end{exm}

\begin{dem}{Proof of Proposition \ref{coarse}.}
Let $\ov{M}_g(\Gamma)$ be the coarse moduli space of $\ov{\M}_g(\Gamma)$, it is a compact algebraic curve. From the universal property of coarse moduli spaces, we have a (non constant) scheme morphism from the coarse moduli space of $\left [\ov{\Ncal}_g(\Gamma) \right]$, that is to say a morphism $\Pro^1 \to \ov{M}_g(\Gamma)$.

We now want to conclude by proving that $\ov{M}_g(\Gamma)$ is smooth. We use the notation of the proof of Proposition \ref{smooth}.

Let $x: \spe (k) \to \ov{M}_g(\Gamma)$ be a geometric point. From properties of coarse moduli spaces, it lifts to a geometric point of $\ov{\M}_g(\Gamma)$. Moreover we can assume that (the lift of) $x$ belongs to $\ov{\M}_g(\Gamma) \backslash \M_g(\Gamma)$. Let $\Gamma_x$\footnote{We get $\Gamma_x$ from a pop of $\Gamma$ in $v_0$.} be the topological class of the curve $C_x$ associated to $x$. The complete local ring of $\ov{M}_g(\Gamma)$ at $x$ is the quotient of $\widehat{\Of}_{x, \ov{\M}_g}={\mathfrak o}_k[[t_1,\ldots , t_{3g-3}]]$ by the ideal corresponding to deformations smoothing a node that allows to get a curve with topological class $\Gamma$ (that is to say a shrink from $\Gamma_x$ to $\Gamma$). After a possible reorder of the variables, we can assume that $t_1, \dots, t_k$ correspond to the nodes $p_1, \dots, p_k$ smoothing which we get a curve of topological class $\Gamma$. Therefore, the ideal we are looking for, is\footnote{This ideal is the intersection of the ideals $I_i:=(t_1, \dots, \widehat{t}_i, \dots , t_{3g-3})$ for $i=1, \dots, k$}
$$
I_x:=\left\{ t_i t_j \right\}_{i<j; i,j=1, \dots, k} \cup \left\{ t_{k+1}, \dots, t_{3g-3} \right\}.
$$
So we get
$$
\widehat{\Of}_{x, \ov{\M}_g(\Gamma)}={\mathfrak o}_k[[t_1,\ldots , t_{k}]]/ \left( \left\{ t_i t_j \right\}_{i<j; i,j=1, \dots, k} \right).
$$

In order to compute the local ring $\widehat{\Of}_{x, \ov{M}_g(\Gamma)}$ of $x$ in the moduli space $\ov{M}_g(\Gamma)$ we have to consider the invariant part of the action of the group $\aut(\Gamma_x)$ (see Proposition \ref{localinvariant}).

Since this action comes from a transitive permutation of the nodes $p_1, \dots, p_k$ (see Lemma \ref{transitive}), through direct computation we get
$$
\widehat{\Of}_{x, \ov{M}_g(\Gamma)} = {\mathfrak o}_k[[s]]
$$
for some generator $s$, and we conclude that $\ov{M}_g(\Gamma)$ is smooth.
\end{dem}

\begin{lemma}\label{transitive}
Let $\Gamma$ a graph of genus $g$ and $3g-4$ edges. Let $\Gamma_x$ be a graph obtained by a pop at the vertex of $\Gamma$ with multiplicity 4. Then $\aut(\Gamma_x)$ acts transitively on the edges $p_1, \dots, p_k$ shrinking which we get the graph $\Gamma$.
\end{lemma}

\begin{dem}
Without loss of generality we show that there exists an element of $\aut(\Gamma_x)$ exchanging $p_1$ and $p_2$.

We know that there exists an isomorphism $\gamma$ between $\Gamma_x$ shrunk at $p_1$ and $\Gamma_x$ shrunk at $p_2$. Clearly $\gamma$ lifts to an automorphism of $\Gamma_x$. 
\end{dem}

\begin{defi}\label{singular}
Let $\Gamma$ be a stable graph of genus $g$ with total weight 0 and $3g-4$ edges. Let 
$$
\spe (k) \xrightarrow{x} \ov{\M}_g(\Gamma) \backslash \M_g(\Gamma)
$$
be a geometric point and $\Gamma_x$ its associated graph. We call order of $x$ in $\ov{\M}_g(\Gamma)$ (written $\ord_x(\Gamma)$) the number of edges in $\Gamma_x$ shrinking which we get the graph $\Gamma$. From the proof of Theorem \ref{coarse} we get that $\ord_x(\Gamma)$ is equal to the number of generators for the complete local ring $\widehat{\Of}_{x, \ov{\M}_g(\Gamma)}$. 
\end{defi}

\begin{rmk} \label{smoothcriterion}
The point $x$ is smooth in $\ov{\M}_g(\Gamma)$ if and only if $ord_x(\Gamma)=1$.
\end{rmk}

\begin{rmk}\label{globalHgerbe}
It is tempting to describe $\ov{\M}_g(\Gamma)$ as an $H$-gerbe over a
Deligne-Mumford stack with only finite stacky points for some finite
group. However, because $\ov{\M}_g(\Gamma)$ is not normal this is not
the case, as we the following example shows.

Let $\Gamma$ be
$$
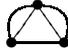
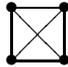
\cthreeone
$$
and $\Gamma_x$ be
$$
\pthreeone
$$
We have $\aut(\Gamma)=D_8$. The group $H$ should be isomorphic to $\sigma^{-1}(V_4)$ (see Theorem \ref{structurenormh0}), which in this case is $V_4 \cong \Z/2\Z \times \Z/2\Z$. Now we have $\aut(\Gamma_x) \cong S_4$, but the group $\sigma^{-1}(V_4)$ seen in $\aut(\Gamma_x)$ is the subgroup fixing an edge which is not normal in $\aut(\Gamma_x)$.
\end{rmk}

\subsection{The case $h(\Gamma)=1$}.

In this case we have $\mathcal N_g(\Gamma)=\mathcal M_{1,1}$. It is known that $\ov{\mathcal M}_{1,1}$ is a $\Z/2\Z$-gerbe over a $[1|2,3]$ orbifold (the coarse moduli space is still $\Pro^1$). Moreover in this case the stabilizer of points is the entire group $\aut(\Gamma)$ and it acts trivially on $\mathcal M_{1,1}$ as it must fix the vertex of weight 1. So we have
$$
[\ov{\mathcal M}_{1,1}/\aut(\Gamma)] = \ov{\mathcal M}_{1,1} \times \B \aut(\Gamma)
$$
and we get 
\begin{theorem}\label{structurenormh1}
The normalization of $\ov{\M}_g(\Gamma)$ is a $(\Z/2\Z \times \aut(\Gamma))$-gerbe over the orbifold $[1|2,3]$.
\end{theorem}

\begin{rmk}\label{caseh1}
Proposition \ref{coarse} still holds and we can also extend Definition \ref{singular}.
\end{rmk}


\subsection{The 1-stratum of $\ov{\Mcal}_3$}

As an example we present the case $g=3$. In the following picture we alternate the graphs with 5 and 6 nodes.

\bigskip

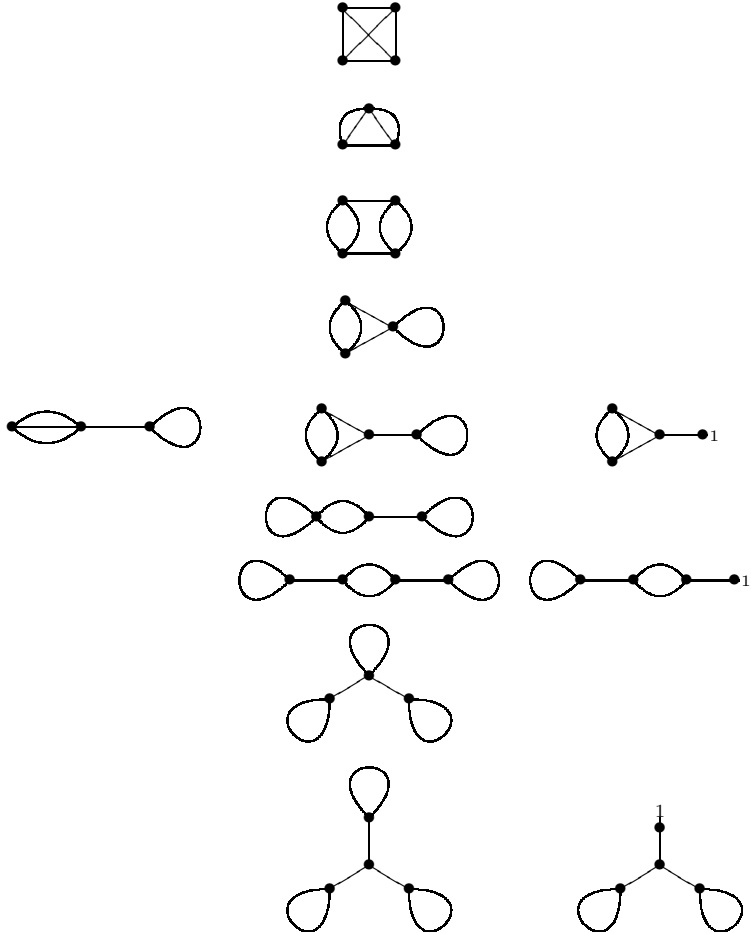
\begin{tabular}{ccc}

& $\pthreeone$ &\\

& &\\

& $\cthreeone$ &\\

&  &\\

& $\pthreetwo$ &\\

&  &\\

& $\cthreetwo$ &\\

&  &\\

$\cthreefive \qquad $ & $\qquad \pthreethree \qquad$ & $\qquad \cthreesix$\\

&  &\\

& $\cthreethree$ & \\

&&\\

& $\qquad \pthreefour \qquad$ & $\qquad \cthreeseven$ \\

&&\\

&&\\

& $\cthreefour$ & \\

&  &\\

&&\\

&&\\

& $\qquad \pthreefive \qquad$  & $\qquad \cthreeeight$ \\

\end{tabular}

\newpage

In the following table we describe the stacks $\ov{\M}_g(\Gamma)$. We use the notation of  \ref{residualorbifolds} and Definition \ref{singular}.

\bigskip

\begin{tabular}{|c|c|c|c|}

\hline
&&& order of points \\
$\Gamma$ & $\sigma^{-1}(V_4)$ &  $\left[ \ov{\N}_{g}(\Gamma)/ \aut(\Gamma) \right]_{\text{orb}}$& in the 0-stratum\\
&&&\\

\hline
$\cthreeone$ & $V_4$ & $\left[ 1, 2 |2 \right]$& $\ord_{x_1}(\Gamma)=6; \ord_{x_2}(\Gamma)=2$\\

\hline
$\cthreetwo \qquad$ & $\left( \Z/2\Z \right)^2$ & $[1,2|2]$ & $\ord_{x_1}(\Gamma)=4; \ord_{x_2}(\Gamma)=1$ \\

\hline

&&&\\

$\qquad \cthreethree \qquad$ & $\left( \Z/2\Z \right)^2$ &  $[1,2|2]$ & $\ord_{x_1}(\Gamma)=2; \ord_{x_2}(\Gamma)=2$\\

&&&\\

\hline

&&&\\

\begin{tabular}{c} $\quad$ \\ $\cthreefour \qquad$ \\ $\quad$ \end{tabular} & $\left( \Z/2\Z \right)^3$ &  $[1,2|2]$ & $\ord_{x_1}(\Gamma)=2; \ord_{x_2}(\Gamma)=3$\\

&&&\\

\hline
$\cthreefive \qquad$ & $\Z/2\Z$ & $[2|2,3]$ & $\ord_{x_1}(\Gamma)=2$\\

\hline

$\qquad \cthreesix$ & $(\Z/2\Z)^3$ & $[1|2,3]$ & $\ord_{x_1}(\Gamma)=1$\\

\hline

&&&\\

$\qquad \cthreeseven$ & $(\Z/2\Z)^3$ & $[1|2,3]$ & $\ord_{x_1}(\Gamma)=2$\\

&&&\\

\hline

\begin{tabular}{c} $\quad$ \\ $\qquad \cthreeeight$ \\ $\quad$ \end{tabular} & $(\Z/2\Z)^4$ &$[1|2,3]$ & $\ord_{x_1}(\Gamma)=3$\\

&&&\\

\hline

\end{tabular}

\begin{exm}\label{allcases}
We want to show that all cases for residual orbifolds exist.

A graph $\Gamma$ such that the residual orbifold of $\left[ \ov{\N}_{g}(\Gamma)/ \aut(\Gamma) \right]$ is $[1,1,1|\emptyset]$, is
$$
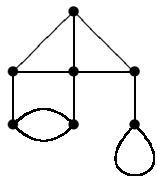
\cfivespecial
$$

A graph $\Gamma$ such that the residual orbifold of $\left[ \ov{\N}_{g}(\Gamma)/ \aut(\Gamma) \right]$ is $[1| 3,3]$, is
$$
\graphonethreethree
$$

\end{exm}


\bigskip














\end{document}